# A UNIVERSAL PROCEDURE FOR AGGREGATING ESTIMATORS

By Alexander Goldenshluger[1]

*University of Haifa*

In this paper we study the aggregation problem that can be formulated as follows. Assume that we have a family of estimators $\mathcal{F}$ built on the basis of available observations. The goal is to construct a new estimator whose risk is as close as possible to that of the best estimator in the family. We propose a general aggregation scheme that is universal in the following sense: it applies for families of arbitrary estimators and a wide variety of models and global risk measures. The procedure is based on comparison of empirical estimates of certain linear functionals with estimates induced by the family $\mathcal{F}$. We derive oracle inequalities and show that they are unimprovable in some sense. Numerical results demonstrate good practical behavior of the procedure.

**1. Introduction.** The subject of this paper is the problem of aggregating estimators from a given collection.

Consider the Gaussian white noise model

$$(1) \quad Y_\varepsilon(dt) = f(t)\,dt + \varepsilon W(dt), \qquad t = (t_1, \ldots, t_d) \in \mathcal{D}_0 = [0,1]^d,$$

where $f : \mathbb{R}^d \to \mathbb{R}$ is an unknown function, $\varepsilon \in (0,1)$ and $W$ is the standard Wiener process in $\mathbb{R}^d$. Let $\Theta \subset \mathbb{R}^N$ be a compact set, and assume that we are given a parameterized family of estimators $\mathcal{F}_\Theta = \{f_\theta, \theta \in \Theta\}$ of $f$. The objective is, using the observation $\mathcal{Y}_\varepsilon = \{Y_\varepsilon(t), t \in \mathcal{D}_0\}$, to select a single estimator from $\mathcal{F}_\Theta$ with the risk that is as close as possible to the risk of the best estimator in the family $\mathcal{F}_\Theta$. We refer to the outlined setup as the *aggregation problem*. Aggregation is a common approach to construction of nonparametric adaptive estimators; this fact motivates consideration of aggregation problems.

Received October 2007; revised November 2007.
[1]Supported in part by ISF Grant 389/07 and by BSF Grant 2006075.
*AMS 2000 subject classifications.* Primary 62G08; secondary 62G05, 62G20.
*Key words and phrases.* Aggregation, lower bound, normal means model, oracle inequalities, sparse vectors, white noise model.

Typically aggregation procedures involve splitting the sample into two sub-samples: the candidate estimators are constructed on the basis of the first sub-sample, while the second subsample is used for the aggregation purposes. In this work we focus on the aggregation step only, and following Juditsky and Nemirovski (2000), Nemirovski (2000) and Tsybakov (2003) we regard the estimators $f_\theta$, $\theta \in \Theta$, as known fixed functions on $\mathcal{D}_0$.

The following two types of aggregation are frequently discussed in the literature:

(i) *Model selection (MS) aggregation.* Here $\Theta = I_N := (1, \ldots, N)$, and the corresponding set of estimators is $\mathcal{F}_\Theta = \mathcal{F}_{I_N} := \{f_i, i \in I_N\}$, where $f_i$ are distinct fixed functions.

(ii) *Convex aggregation.* Here

$$\Theta = \Lambda := \left\{ \lambda \in \mathbb{R}^N | \lambda_i \geq 0, \sum_{i=1}^N \lambda_i \leq 1 \right\}, \tag{2}$$

and for fixed estimators $f_i$, $i \in I_N$,

$$\mathcal{F}_\Theta = \mathcal{F}_\Lambda := \left\{ F_\lambda | F_\lambda(t) := \sum_{i=1}^N \lambda_i f_i(t), \lambda \in \Lambda \right\}.$$

Let $\tilde{f}$ be an estimator of $f$ based on the observation $\mathcal{Y}_\varepsilon$. We measure accuracy of $\tilde{f}$ by its $\mathbb{L}_p$-risk

$$\mathcal{R}_p[\tilde{f}; f] := \mathbb{E}_f \|\tilde{f} - f\|_p, \qquad 1 \leq p \leq \infty,$$

where $\mathbb{E}_f$ is the expectation with respect to the probability measure $\mathbb{P}_f$ of observation $\mathcal{Y}_\varepsilon$ under model (1), and $\|\cdot\|_p$ is the standard $\mathbb{L}_p$-norm on $\mathcal{D}_0$. We want to propose a measurable choice, say $\hat{f} = f_{\hat{\theta}}$, from collection $\mathcal{F}_\Theta$ such that the following $\mathbb{L}_p$-*risk oracle inequality* holds:

$$\mathcal{R}_p[\hat{f}; f] \leq C \inf_{\theta \in \Theta} \mathcal{R}_p[f_\theta; f] + r_\varepsilon \tag{3}$$

for all $f$ from a "large" functional class. Here $C$ is a constant independent of $f$ and $\varepsilon$, and $r_\varepsilon$ is a remainder term that does not depend on $f$.

The outlined aggregation problem has attracted much attention in the literature for the regression and Gaussian white noise models. Remarkable progress has been achieved in the framework of $\mathbb{L}_2$-theory where exact oracle inequalities [with $C = 1$ or $C = 1 + o(1)$, $\varepsilon \to 0$] were derived for collections of arbitrary estimators; see Juditsky and Nemirovski (2000), Nemirovski (2000), Tsybakov (2003). Tsybakov (2003) introduced the notion of optimal rates of aggregation and derived aggregation procedures possessing (3) with smallest possible, in a minimax sense, remainder term $r_\varepsilon$. $\mathbb{L}_2$-risk oracle inequalities with $C > 1$ for arbitrary estimators were obtained, for example, by Yang (2001, 2004), Wegkamp (2003) and Bunea, Tsybakov and Wegkamp (2007).

Aggregation of arbitrary nonparametric estimators with respect to other loss functions is much less studied. Catoni (2004) and Yang (2000) considered the problem of aggregating density estimators with the Kullback–Leibler divergence as a loss function. Devroye and Lugosi (1996, 1997, 2001) developed $\mathbb{L}_1$-risk oracle inequalities in the context of density estimation; see also Hengartner and Wegkamp (2001) who apply the approach of Devroye and Lugosi for the regression setup. Our results are closely related to those by Devroye and Lugosi, and we discuss this connection in detail in Section 3.

For a detailed account of the literature on aggregation of estimators see the recent papers Audibert (2004), Birgé (2006), Bunea, Tsybakov and Wegkamp (2007), Juditsky, Rigollet and Tsybakov (2008) and references therein. It is also worth noting that there is vast literature on aggregation of estimators from restricted families (such as orthogonal series estimators, kernel estimators, etc.), and aggregation of classifiers in classification problems. A list of representative publications from this literature includes Kneip (1994), Lepski and Spokoiny (1997), Cavalier et al. (2002), Koltchinskii (2006) and Lecué (2007), where further references can be found.

In this paper we propose a general aggregation scheme that is universal in the following sense: (i) it applies to families of arbitrary estimators; (ii) it can be easily extended to different models; (iii) it can be used for a wide variety of global risk measures. Although the main results of this paper pertain to the MS aggregation setup, Gaussian white noise model and $\mathbb{L}_p$-risks, similar results can be easily established for other models and global risk measures. In Section 4 we illustrate universality of the suggested procedure by applying it to convex aggregation and to the problem of estimating a normal mean vector.

Our aggregation method is based on comparison of empirical estimates of certain regular linear functionals with estimates induced by the family $\mathcal{F}_\Theta$. A closely related idea that a nonparametric function estimator is "good" if its integrals over cubes "agree" with the corresponding empirical means, belongs to Nemirovski (1985). We establish general oracle inequalities and specialize them for different sets of linear functionals. It turns out that universal inequalities of Devroye and Lugosi (1996, 1997, 2001) and Hengartner and Wegkamp (2001) can be derived from our general oracle inequalities using a specific choice of the set of linear functionals. The results indicate that in the Gaussian white noise model (1) the problem of aggregation of *arbitrary* estimators in $\mathbb{L}_p$, $p \in (2, \infty]$, can be rather difficult. In this case remainder terms in the oracle inequalities depend on the family $\mathcal{F}_\Theta$ and, in general, can be rather large. We prove a lower bound and show that dependence of the remainder terms on $\mathcal{F}_\Theta$ is, in a sense, unavoidable.

Thus "efficient" aggregation of *arbitrary* estimators in $\mathbb{L}_p$, $p \in (2, \infty]$, is impossible. We also show that in the $\mathbb{L}_2$-framework a slight modification of the proposed aggregation procedure satisfies the exact oracle inequality (3) with $C = 1$ and the remainder $r_\varepsilon$ that cannot be improved in the minimax sense.

The rest of the paper is organized as follows. In Section 2 we introduce our aggregation scheme. Section 3 contains the main results of the paper. In Section 4 we apply the procedure to convex aggregation and estimation of a normal mean vector. In a simulation experiment of Section 4 we study performance of our procedure for estimating a normal mean vector. Proofs are given in Section 5.

**2. Aggregation scheme.** We begin with construction of the aggregation scheme for the Gaussian white noise model (1).

2.1. *Construction.* Let $\Psi$ be a set of probe functions $\psi : \mathcal{D}_0 \to \mathbb{R}$. Consider a linear functional
$$\ell_f(\psi) = \int \psi(t) f(t) \, dt, \qquad \psi \in \Psi.$$

For given $\psi \in \Psi$, a natural estimator of $\ell_f(\psi)$ based on observation $\mathcal{Y}_\varepsilon$ is
$$\hat{\ell}_f(\psi) = \int \psi(t) Y_\varepsilon(dt).$$

On the other hand, $\ell_f(\psi)$ can be estimated using estimates $f_\theta \in \mathcal{F}_\Theta$:
$$\ell_{f_\theta}(\psi) = \int \psi(t) f_\theta(t) \, dt, \qquad \theta \in \Theta.$$

Define
$$\begin{aligned}
\Delta_\theta(\psi) &:= \hat{\ell}_f(\psi) - \ell_{f_\theta}(\psi) \\
&= \int \psi(t)[f(t) - f_\theta(t)] \, dt + \varepsilon \int \psi(t) W(dt) \\
(4) \qquad &=: \int \psi(t)[f(t) - f_\theta(t)] \, dt + \varepsilon Z(\psi), \qquad \theta \in \Theta.
\end{aligned}$$

For any fixed $\theta \in \Theta$, $\Delta_\theta(\psi)$ is a random variable that measures discrepancy between empirical estimate $\hat{\ell}_f(\psi)$ of the linear functional $\ell_f(\psi)$ and the estimate $\ell_{f_\theta}(\psi)$ induced by $f_\theta \in \mathcal{F}_\Theta$. The idea underlying construction of our aggregation rule is that, for a "good" estimator $f_\theta$, the absolute value of $\Delta_\theta(\psi)$ "corrected" for a random error $Z(\psi)$ should be uniformly small for all $\psi \in \Psi$.

Let $\delta \in (0, 1)$, and
$$(5) \qquad \varkappa = \varkappa(\delta, \Psi) := \min\left\{ \varkappa > 0 \Big| \mathbb{P}\left[\sup_{\psi \in \Psi} \frac{|Z(\psi)|}{\|\psi\|_2} \geq \varkappa\right] \leq \delta \right\}.$$

Define

$$\hat{M}_\theta := \sup_{\psi \in \Psi} \left\{ \frac{1}{\|\psi\|_q} [|\Delta_\theta(\psi)| - \varepsilon \varkappa \|\psi\|_2] \right\}, \tag{6}$$

where $p^{-1} + q^{-1} = 1$, and let $\hat{\theta} := \arg\inf_{\theta \in \Theta} \hat{M}_\theta$; then our estimator is given by

$$\hat{f} = f_{\hat{\theta}}. \tag{7}$$

Recently a procedure based on different ideas but close in spirit to (6)–(7) was used in Goldenshluger and Lepski (2007) for selection of kernel estimators from large parameterized collections.

In order to ensure that the estimator $\hat{f}$ is well defined, certain conditions on the set of probe functions $\Psi$, and on the family of estimators $\mathcal{F}_\Theta$, have to be imposed. First, to guarantee that $\varkappa$ is well defined in (5), we need appropriate assumptions on the intrinsic semimetric of the zero-mean Gaussian process $\{Z(\psi), \psi \in \Psi\}$. Second, $\hat{\theta}$ should be measurable; this requirement calls for conditions on the sample paths of the random process $\{\hat{M}_\theta, \theta \in \Theta\}$. Although general conditions that guarantee fulfillment of the above properties can be explicitly stated, for the present we will take them for granted. In the aggregation setups of Sections 3 and 4 these conditions are trivially fulfilled.

Note that the aggregation procedure requires specification of the parameter $\delta$ and the set of probe functions $\Psi$. The choice of $\Psi$ is a crucial step in construction. We discuss this issue below.

2.2. *The set of probe functions.* The following *norm approximation property* of the set of probe functions $\Psi$ plays an important role in our construction.

DEFINITION 1. Given the collection of estimators $\mathcal{F}_\Theta = \{f_\theta, \theta \in \Theta\}$ with index set $\Theta$, let

$$\mathcal{G}_\Theta := \{g : \mathcal{D}_0 \to \mathbb{R} | g = g_{\tau,\nu} := f_\tau - f_\nu, f_\tau, f_\nu \in \mathcal{F}_\Theta, f_\tau \neq f_\nu\}. \tag{8}$$

Let $\Psi$ be a set of functions on $\mathcal{D}_0$, $\gamma \geq 0$ and $p \in [1, \infty]$. We say that $\Psi$ is a $(\gamma, p)$-*good set with respect to* $\mathcal{G}_\Theta$ if for any $g \in \mathcal{G}_\Theta$ there exists $\psi_g \in \Psi$ such that

$$\left| \int \psi_g(t) g(t) \, dt - \|g\|_p \right| \leq \gamma. \tag{9}$$

Several remarks on the above definition are in order. The set $\mathcal{G}_\Theta$ contains pairwise differences of estimators from $\mathcal{F}_\Theta$. The set of probe functions $\Psi$ is $(\gamma, p)$-good with respect to $\mathcal{G}_\Theta$ if the $\mathbb{L}_p$-norm of any function from $\mathcal{G}_\Theta$ can be approximated by a linear functional from $\Psi$ with prescribed guaranteed accuracy $\gamma$. Since $\mathcal{G}_\Theta$ is indexed by $(\tau,\nu) \in \Theta \times \Theta$, the corresponding $(\gamma,p)$-good set of probe functions can be always chosen indexed by $(\tau,\nu) \in \Theta \times \Theta$, too. Specifically, the $(\gamma,p)$-good set with respect to $\mathcal{G}_\Theta$ can be chosen as follows:

$$\Psi = \Psi_\Theta := \{\psi : \mathcal{D}_0 \to \mathbb{R} | \psi = \psi_{g_{\tau,\nu}}, \tau, \nu \in \Theta,\ \tau \neq \nu\}, \tag{10}$$

where $\psi_{g_{\tau,\nu}}$ is the representer corresponding to $g_{\tau,\nu} := f_\tau - f_\nu$ such that (9) is fulfilled. In all that follows without further mention we always write $\Psi_\Theta$ for a set of probe functions that is associated with $\Theta$ (and $\mathcal{G}_\Theta$) via (10).

The $(\gamma,p)$-good sets of probe functions are easily constructed. In the sequel the following examples of the $(\gamma,p)$-good sets will be particularly important.

EXAMPLE 1. Let $p \in [1,\infty)$ and define

$$\widetilde{\Psi}_\Theta := \left\{\psi | \psi(\cdot) = \psi_g(\cdot) := \frac{|g(\cdot)|^{p-1}}{\|g\|_p^{p-1}} \operatorname{sign}\{g(\cdot)\}, g \in \mathcal{G}_\Theta\right\}. \tag{11}$$

Clearly, $\widetilde{\Psi}_\Theta$ is $(0,p)$-good with respect to $\mathcal{G}_\Theta$. Note also that $\widetilde{\Psi}_\Theta \subseteq \{\psi : \|\psi\|_q = 1\}$.

EXAMPLE 2. The set

$$\widehat{\Psi}_\Theta := \left\{\psi | \psi(\cdot) = \psi_g(\cdot) := \frac{\|g\|_p}{\|g\|_2^2} g(\cdot), g \in \mathcal{G}_\Theta\right\} \tag{12}$$

is $(0,p)$-good with respect to $\mathcal{G}_\Theta$ for any $p \in [1,\infty]$.

EXAMPLE 3. For $\gamma > 0$ define

$$\overline{\Psi}_\Theta(\gamma) := \left\{\psi | \psi(\cdot) = \psi_g(\cdot) := \frac{[|g(\cdot)| - \|g\|_\infty + \gamma]_+ \operatorname{sign}\{g(\cdot)\}}{\int [|g(t)| - \|g\|_\infty + \gamma]_+\, dt}, g \in \mathcal{G}_\Theta\right\},$$

where $[\cdot]_+ = \max\{\cdot, 0\}$. It is easily verified that $\overline{\Psi}_\Theta(\gamma)$ is $(\gamma,\infty)$-good with respect to $\mathcal{G}_\Theta$; moreover, $\overline{\Psi}_\Theta(\gamma) \subset \{\psi : \|\psi\|_1 = 1\}$.

**3. Main results.** In this section we present the main results of this paper. We focus on the model selection aggregation setup where $\Theta = I_N = (1,\ldots,N)$, $\mathcal{F}_\Theta = \mathcal{F}_{I_N} = \{f_i, i \in I_N\}$. Let $\mathcal{G}_{I_N}$ and $\Psi_{I_N}$ be defined accordingly via (8) and (10). Note that $\mathcal{G}_{I_N}$ and $\Psi_{I_N}$ are finite sets of functions of cardinality $N(N-1)$. Following (4), for $\psi \in \Psi_{I_N}$ we write

$$\Delta_i(\psi) := \hat{\ell}_f(\psi) - \ell_{f_i}(\psi)$$

$$= \int \psi(t)[f(t) - f_i(t)]\,dt + \varepsilon Z(\psi), \qquad i \in I_N. \tag{13}$$

For a fixed $\delta \in (0,1)$, $\varkappa = \varkappa(\delta, \Psi_{I_N})$ is given by (5); note that $\varkappa$ is well defined because $\Psi_{I_N}$ is a finite set. We write also

$$\hat{M}_i := \max_{\psi \in \Psi_{I_N}} \left\{ \frac{1}{\|\psi\|_q} [|\Delta_i(\psi)| - \varepsilon\varkappa\|\psi\|_2] \right\} \tag{14}$$

and

$$\hat{i} := \arg\min_{i \in I_N} \hat{M}_i, \qquad \hat{f} = f_{\hat{i}}. \tag{15}$$

3.1. *Oracle inequalities.* The next theorem establishes the basic oracle inequality on the $\mathbb{L}_p$-risk of the estimator $\hat{f}$.

THEOREM 1. *Let $p \in [1,\infty]$, and assume that $\Psi_{I_N}$ is $(\gamma,p)$-good with respect to $\mathcal{G}_{I_N}$. Define $i_* := \arg\min_{i \in I_N} \|f - f_i\|_p$ and*

$$\Psi^*_{I_N} := \{\psi \in \Psi_{I_N} | \psi = \psi_{f_{i_*} - f_i} = \psi_{i_* i}, i \in I_N, i \neq i_*\}. \tag{16}$$

*Let $\delta \in (0,1)$ be fixed, and let $\varkappa = \varkappa(\delta, \Psi_{I_N})$ be defined in (5); then for $\hat{f}$ given in (14)–(15) one has*

$$\begin{aligned}
\mathcal{R}_p[\hat{f}; f] &\leq \left(2 \max_{\psi \in \Psi^*_{I_N}} \|\psi\|_q + 1\right) \min_{i \in I_N} \|f - f_i\|_p \\
&\quad + 2\varkappa\varepsilon \max_{\psi \in \Psi^*_{I_N}} \|\psi\|_2 + \gamma + \left[\|f\|_p + \max_{i \in I_N} \|f_i\|_p\right]\delta.
\end{aligned} \tag{17}$$

REMARK 1. The proof of Theorem 1 illuminates the role played by the assumption that $\Psi_{I_N}$ is $(\gamma,p)$-good. The key is the bound on the distance between selected and oracle estimators, $\|f_{i_*} - f_{\hat{i}}\|_p$. The fact that $\Psi_{I_N}$ is $(\gamma,p)$-good allows to control this distance on an event of large probability in terms of the distance between corresponding linear functionals. The latter, in turn, is controlled by definition of the aggregation procedure.

We now apply the oracle inequality of Theorem 1 for the sets of probe functions discussed in Examples 1–3 of Section 2. Assume that

$$\max\{\|f\|_p, \|f_1\|_p, \ldots, \|f_N\|_p\} := L < \infty. \tag{18}$$

COROLLARY 1. *Let $\Psi_{I_N} = \widetilde{\Psi}_{I_N}$ where $\widetilde{\Psi}_\Theta$ is defined in (11). Suppose that (18) holds; then for $\hat{f}$ given in (14)–(15) and associated with $\widetilde{\Psi}_{I_N}$ and $\delta = \varepsilon$ one has*

$$\mathcal{R}_p[\hat{f}; f] \leq 3 \min_{i \in I_N} \|f - f_i\|_p + 2Q_1(p)\varepsilon\sqrt{2\ln\frac{N^2}{\varepsilon}} + 2L\varepsilon, \tag{19}$$

where $Q_1(p) = 1$ for $1 \leq p \leq 2$, and

$$(20) \quad Q_1(p) = Q_1(\mathcal{F}_{I_N}, p) := \max_{\substack{i \in I_N \\ i \neq i_*}} \left[ \frac{\|f_{i_*} - f_i\|_{2p-2}}{\|f_{i_*} - f_i\|_p} \right]^{p-1}, \quad 2 < p < \infty.$$

REMARK 2. Our selection rule with $\Psi_{I_N} = \widetilde{\Psi}_{I_N}$ and $p = 1$ reduces to the aggregation method by Devroye and Lugosi (1996, 1997, 2001). Indeed, when $p = 1$, the probe functions from the set $\widetilde{\Psi}_{I_N}$ are given by $\psi_{ij} = \text{sign}(f_i - f_j)$. In the density estimation context this corresponds to the Yatracos classes considered by Devroye and Lugosi. Note also that when $p \in [1, 2]$ and $\Psi_{I_N} = \widetilde{\Psi}_{I_N}$, the selection rule (14)–(15) could be modified as follows:

$$\hat{i} = \arg\min_{i \in I_N} \max_{\psi \in \widetilde{\Psi}_{I_N}} |\Delta_i(\psi)|.$$

In this form our selection rule can be viewed as an implementation of the method by Devroye and Lugosi for the white noise model [see also Hengartner and Wegkamp (2001)]. For further discussion see Section 3.3.

COROLLARY 2. Let $p \in [1, \infty]$, and $\Psi = \widehat{\Psi}_{I_N}$ where $\widehat{\Psi}_\Theta$ is defined in (12). Suppose that (18) holds; then for the estimate $\hat{f}$ given in (14)–(15) and associated with $\widehat{\Psi}_{I_N}$ and $\delta = \varepsilon$ one has

$$(21) \quad \mathcal{R}_p[\hat{f}; f] \leq (2Q_2(p) + 1) \min_{i \in I_N} \|f_i - f\|_p + 2Q_3(p)\varepsilon \sqrt{2\ln \frac{N^2}{\varepsilon}} + 2L\varepsilon,$$

where

$$(22) \quad \begin{aligned} Q_2(p) &= Q_2(\mathcal{F}_{I_N}, p) := \max_{\substack{i \in I_N \\ i \neq i_*}} \frac{\|f_{i_*} - f_i\|_p \|f_{i_*} - f_i\|_q}{\|f_{i_*} - f_i\|_2^2}, \\ Q_3(p) &= Q_3(\mathcal{F}_{I_N}, p) := \max_{\substack{i \in I_N \\ i \neq i_*}} \frac{\|f_{i_*} - f_i\|_p}{\|f_{i_*} - f_i\|_2}. \end{aligned}$$

In contrast to $\widetilde{\Psi}_{I_N}$, the rule associated with $\widehat{\Psi}_{I_N}$ allows to treat the case $p = \infty$. Note, however, that it leads to the elevated factor preceding the best possible risk as compared to the selection rule that uses $\widetilde{\Psi}_{I_N}$.

COROLLARY 3. Let (18) hold with $p = \infty$, and $\Psi_{I_N} = \overline{\Psi}_{I_N}(\gamma_0)$ with $\gamma_0 = \varepsilon\sqrt{\ln N} < L$; then

$$(23) \quad \mathcal{R}_\infty[\hat{f}; f] \leq 3 \min_{i \in I_N} \|f_i - f\|_\infty + 3Q_4(\gamma_0)\varepsilon \sqrt{2\ln \frac{N^2}{\varepsilon}} + 2L\varepsilon,$$

*where*

$$Q_4(\gamma) = Q_4(\mathcal{F}_{I_N}, \gamma) := \max_{\substack{i \in I_N \\ i \neq i_*}} \frac{\|S_{i_*i}(\cdot, \gamma)\|_2}{\|S_{i_*i}(\cdot, \gamma)\|_1},$$

(24)
$$S_{i_*i}(\cdot, \gamma) := [|f_{i_*}(\cdot) - f_i(\cdot)| - \|f_{i_*} - f_i\|_\infty + \gamma]_+.$$

The above results show that when $p \in [1, 2]$ arbitrary estimators satisfying (18) can be *efficiently* aggregated in the following sense. Corollary 1 demonstrates that if $\Psi = \widetilde{\Psi}_{I_N}$, then the resulting risk of the selected estimator is within factor 3 of the best possible risk whereas the remainder term is of the order $\varepsilon\sqrt{\ln(N^2/\varepsilon)}$. Thus one can aggregate polynomial in $\varepsilon^{-1}$ number $N$ of estimators with remainder term of the order $\varepsilon\sqrt{\ln(1/\varepsilon)}$. Such a bound allows to derive minimax and adaptive results in many nonparametric estimation setups.

The situation is completely different for $p \in (2, \infty]$. Here remainder terms in the oracle inequalities depend on the family of aggregated estimates through the values of $Q_1(p)$, $Q_3(p)$ and $Q_4(\gamma)$ that can be large for particular families $\mathcal{F}_{I_N}$.

3.2. *Lower bound.* The important question is whether the remainder terms in (19), (21) and (23) can be improved for families of arbitrary estimators $\mathcal{F}_{I_N}$ whenever $p > 2$. The next result shows that, in a sense, dependence of the remainder terms on the family $\mathcal{F}_{I_N}$ is unimprovable in the MS aggregation setup.

THEOREM 2. *Assume that $N > 3$ and $p \in (2, \infty]$; then there exists a family $\bar{\mathcal{F}}_{I_N} = \{\bar{f}_i, i \in I_N\}$ of functions on $\mathcal{D}_0$, satisfying $\max_{i \in I_N} \|\bar{f}_i\|_p \leq L$ such that for any selection rule $\tilde{f} : \mathcal{Y}_\varepsilon \to \bar{\mathcal{F}}_{I_N}$ and any $\varepsilon \leq L(N \ln N)^{-1/2}$ one has*

(25) $$\max_{f \in \bar{\mathcal{F}}_{I_N}} \left[ \mathcal{R}_p[\tilde{f}; f] - \min_{i \in I_N} \|f - \bar{f}_i\|_p \right] \geq cK_p \varepsilon \sqrt{\ln(N-1)},$$

*where* $K_p = Q_1(\bar{\mathcal{F}}_{I_N}, p) = Q_3(\bar{\mathcal{F}}_{I_N}, p)$, $\forall p \in [2, \infty)$, $K_\infty = Q_3(\bar{\mathcal{F}}_{I_N}, \infty) = Q_4(\bar{\mathcal{F}}_{I_N}, \gamma)$, $\forall \gamma > 0$, *and $c$ is an absolute constant. The quantities $Q_1$, $Q_3$ and $Q_4$ are defined in (20), (22) and (24), respectively.*

REMARK 3. Because $\min_{i \in I_N} \|f - \bar{f}_i\|_p = 0$ for $f \in \bar{\mathcal{F}}_{I_N}$, (25) provides a lower bound on the remainder term in the $\mathbb{L}_p$-risk oracle inequality. The worst-case family $\bar{\mathcal{F}}_{I_N}$ in the proof of Theorem 2 is such that the $\mathbb{L}_2$-norm of pairwise differences of its members is small in comparison with their $\mathbb{L}_p$-norm. We note also that the worst-case family $\bar{\mathcal{F}}_{I_N}$ does not depend on $p$.

Theorem 2 shows that the problem of aggregation of *arbitrary* estimators in $\mathbb{L}_p$, $p \in (2, \infty]$ may be rather difficult. In particular, the proof of the theorem suggests that the $\mathbb{L}_p$-risk of any aggregation procedure can be as large as $\varepsilon^{2/p}(\ln N)^{1/p}$, $p \in (2, \infty]$.

The meaning of the lower bound of Theorem 2 is that there is a family of estimators that cannot be aggregated with accuracy better than that in (25). This, however, does not imply that the same lower bound holds for a concrete family of reasonable estimators. It is known, for example, that *kernel estimators* can be efficiently aggregated in $\mathbb{L}_p$, $p > 2$ [Goldenshluger and Lepski (2007)].

3.3. *Modified aggregation procedure.* In the definition of the aggregation procedure [see (14)], the "typical" value of the stochastic error, $\varepsilon \varkappa \|\psi\|_2$, is subtracted from $|\Delta_i(\psi)|$. Thus, this construction requires prior knowledge of the noise level $\varepsilon$. We note, however, that the original procedure can be modified in such a way that $\varepsilon$ need not be known.

Specifically, consider the following procedure: with $\Delta_i(\psi)$ given in (13) define

$$\tilde{M}_i := \max_{\psi \in \Psi_{I_N}} \left\{ \frac{1}{\|\psi\|_q} |\Delta_i(\psi)| \right\} \tag{26}$$

and let

$$\tilde{i} := \arg\min_{i \in I_N} \tilde{M}_i, \qquad \tilde{f} = f_{\tilde{i}}. \tag{27}$$

This construction does not require prior knowledge of the noise level $\varepsilon$. The next theorem establishes an oracle inequality for the estimator $\tilde{f}$.

THEOREM 3. *Let conditions of Theorem 1 hold; then for the estimator $\tilde{f}$ defined in (26)–(27) one has*

$$\begin{aligned}
\mathcal{R}_p[\tilde{f}; f] &\leq \left(2 \max_{\psi \in \Psi_{I_N}^*} \|\psi\|_q + 1\right) \min_{i \in I_N} \|f - f_i\|_p \\
&\quad + 2\varkappa\varepsilon \max_{\psi \in \Psi_{I_N}^*} \|\psi\|_q \max_{\psi \in \Psi_{I_N}} \{\|\psi\|_2/\|\psi\|_q\} \\
&\quad + \gamma + \left[\|f\|_p + \max_{i \in I_N} \|f_i\|_p\right]\delta.
\end{aligned} \tag{28}$$

REMARK 4. The second term on the right-hand side of (28) is greater than or equal to that on the right-hand side of (17). However, in special cases oracle inequality (28) is precise enough. For instance, if $p = 2$, then the remainder terms in (28) and (17) coincide. Note also that in the setup of Devroye and Lugosi (2001) ($p = 1$ and $\Psi_{I_N} = \widetilde{\Psi}_{I_N}$; see Remark 2) we obtain

$$2\varkappa\varepsilon \max_{\psi \in \widetilde{\Psi}_{I_N}^*} \|\psi\|_\infty \max_{\psi \in \widetilde{\Psi}_{I_N}} \frac{\|\psi\|_2}{\|\psi\|_\infty} \leq 2\varkappa\varepsilon$$

because $\|\psi\|_2 \leq \|\psi\|_\infty = 1$ for every $\psi \in \widetilde{\Psi}_{I_N}$ whenever $p = 1$. In these cases the use of the modified selection rule is advantageous as it does not require knowledge of the noise level $\varepsilon$.

3.4. $\mathbb{L}_2$-*risk oracle inequality.* If $p = 2$, then the general oracle inequality of Theorem 1 can be improved. In particular, we demonstrate that in this specific case a mild modification of the original aggregation procedure leads to the exact oracle inequality with the leading constant equal to 1.

First we note that the sets of probe functions $\widetilde{\Psi}_{I_N}$ and $\widehat{\Psi}_{I_N}$ coincide when $p = 2$:

$$\psi_{ij}(\cdot) = \frac{f_i(\cdot) - f_j(\cdot)}{\|f_i - f_j\|_2}, \qquad i, j \in I_N, \ i \neq j. \tag{29}$$

Let $u_{ij} = \frac{1}{2}(f_i + f_j)$, and for all $i \in I_N$ define

$$\begin{aligned}
\bar{M}_i &:= \max_{j \in I_N} \{\ell_{u_{ij}}(\psi_{ij}) - \hat{\ell}_f(\psi_{ij})\} \\
&= \max_{j \in I_N} \left\{ \int \psi_{ij}(t) u_{ij}(t)\, dt - \int \psi_{ij}(t) Y_\varepsilon(dt) \right\}.
\end{aligned} \tag{30}$$

The selection rule is defined by

$$\bar{i} = \arg\min_{i \in I_N} \bar{M}_i, \qquad \bar{f} = f_{\bar{i}}. \tag{31}$$

We remark that $\|\psi_{ij}\|_2 = 1$, $\forall i, j \in I_N$, $i \neq j$. A distinctive feature of the selection rule (29)–(31) is that for each pair $i, j \in I_N$ the empirical estimate of the linear functional $\ell_f(\psi_{ij})$ is compared with $\ell_{u_{ij}}(\psi_{ij})$ and not with $\ell_{f_i}(\psi_{ij})$ as in (13).

THEOREM 4. *Let $\bar{f} = f_{\bar{i}}$ be the estimator defined by (29)–(31); then*

$$\mathcal{R}_2[\bar{f}; f] \leq \min_{i \in I_N} \|f_i - f\|_2 + 8\varepsilon\sqrt{2 \ln N}.$$

Thus the selection rule (29)–(31) achieves the optimal rates of the MS aggregation when the $\mathbb{L}_2$-risk is considered [cf. Tsybakov (2003)].

**4. Miscellaneous extensions and numerical results.** The objective of this section is to demonstrate that the proposed procedure can be applied for different models and global risk measures. First we discuss the problem of convex aggregation, and then we show how the aggregation scheme can be applied for estimation in the normal means model. We also provide some numerical results for the problem of estimating a normal mean vector.

4.1. *Convex aggregation.* The problem of convex aggregation is formulated as follows: given a set of estimators $f_i$, $i \in I_N$, the objective is to select an estimator, say $\hat{F} = F_{\hat{\lambda}}$, from the collection

$$\mathcal{F}_\Lambda = \left\{ F_\lambda | F_\lambda(t) = \sum_{i=1}^N \lambda_i f_i(t), \lambda \in \Lambda \right\},$$

such that $F_{\hat{\lambda}}$ is nearly as good as the best estimator from $\mathcal{F}_\lambda$. Here $\Lambda$ is the $N$-dimensional simplex; see (2).

For $\eta > 0$ let $\Lambda_\eta = (\lambda^{(k)}, k = 1, \ldots, n_\eta)$ denote the minimal $\eta$-net of $\Lambda$ in $l_1$-norm; that is, for any $\lambda \in \Lambda$ there exists $\lambda^{(k)} \in \Lambda_\eta$ such that

$$|\lambda - \lambda^{(k)}|_1 = \sum_{i=1}^N |\lambda_i - \lambda_i^{(k)}| \leq \eta.$$

Let $\mathcal{G}_\Lambda = \{g | g = F_\lambda - F_\nu, \lambda, \nu \in \Lambda, \nu \neq \lambda\}$, and let $\mathcal{G}_{\Lambda_\eta}$ be defined similarly with $\Lambda$ replaced by $\Lambda_\eta$ [cf. (8)]. Note that $\mathcal{G}_{\Lambda_\eta}$ is a finite set with $\text{card}(\mathcal{G}_{\Lambda_\eta}) = n_\eta(n_\eta - 1)$.

We begin with a lemma showing that if (18) holds, then any $(0, p)$-good set with respect to $\mathcal{G}_{\Lambda_\eta}$ is also $(\gamma, p)$-good with respect to $\mathcal{G}_\Lambda$ with some $\gamma = \gamma(\eta) > 0$.

LEMMA 1. *Assume that (18) holds, and let $\Psi$ be the $(0, p)$-good set with respect to $\mathcal{G}_{\Lambda_\eta}$. Then $\Psi$ is $(\gamma, p)$-good with respect to $\mathcal{G}_\Lambda$ with*

(32) $$\gamma = 2L\eta \left( 1 + \max_{\psi \in \Psi} \|\psi\|_q \right).$$

Lemma 1 allows to reduce the problem of convex aggregation to the MS aggregation over a finite family of estimators. The idea is to apply the selection procedure of Section 2 to the finite set of estimators induced by the minimal $\eta$-net $\Lambda_\eta$ in $\Lambda$.

Similarly to (13), for $\psi \in \Psi$ we write

$$\Delta_\lambda(\psi) = \hat{\ell}_f(\psi) - \ell_{F_\lambda}(\psi)$$
$$= \int \psi(t)[f(t) - F_\lambda(t)] \, dt + \varepsilon \int \psi(t) W(dt), \qquad \lambda \in \Lambda.$$

Let $\eta = \varepsilon$, and $\Lambda_\varepsilon = \{\lambda^{(k)}, k = 1, \ldots, n_\varepsilon\}$ be a minimal $\varepsilon$-net in $l_1$-norm for $\Lambda$. Let $\Psi_{\Lambda_\varepsilon}$ be a $(0, p)$-good set w.r.t. $\mathcal{G}_{\Lambda_\varepsilon}$. For $\delta \in (0, 1)$ let $\varkappa = \varkappa(\delta, \Psi_{\Lambda_\varepsilon})$ be given by (5). Define

$$(33) \quad \hat{\lambda} := \arg\min_{\lambda \in \Lambda} \max_{\psi \in \Psi_{\Lambda_\varepsilon}} \left\{ \frac{1}{\|\psi\|_2} [|\Delta_\lambda(\psi)| - \varepsilon \varkappa \|\psi\|_2] \right\}, \qquad \hat{F} := F_{\hat{\lambda}}.$$

THEOREM 5. *Assume that $\Psi_{\Lambda_\varepsilon}$ is $(0, \gamma)$-good with respect to $\mathcal{G}_{\Lambda_\varepsilon}$. Then for $\varkappa = \varkappa(\delta, \Psi_{\Lambda_\varepsilon})$ defined in (5) and $\hat{F}$ given by (33) one has*

$$\mathcal{R}_p[\hat{F}; f] \leq \left( 2 \max_{\psi \in \Psi_{\Lambda_\varepsilon}} \|\psi\|_q + 1 \right) \min_{\lambda \in \Lambda} \|f - F_\lambda\|_p$$
$$+ 2\varkappa \varepsilon \max_{\psi \in \Psi_{\Lambda_\varepsilon}} \|\psi\|_2 + 2L\varepsilon \left( 1 + \max_{\psi \in \Psi_{\Lambda_\varepsilon}} \|\psi\|_q \right) + 2L\delta.$$

The oracle inequality of Theorem 5 can be straightforwardly specialized for specific sets of probe functions. We provide here only one result corresponding to Example 1 in Section 2.

COROLLARY 4. *Let $\Psi_{\Lambda_\varepsilon} = \widetilde{\Psi}_{\Lambda_\varepsilon}$ where $\widetilde{\Psi}_\Theta$ is defined in (11). Then for the estimator $\hat{F}$ associated with $\delta = \varepsilon$ one has*

$$\mathcal{R}_p[\hat{F}; f] \leq 3 \min_{\lambda \in \Lambda} \|f - F_\lambda\|_p + cQ_1(p)\varepsilon\sqrt{N \ln \frac{1}{\varepsilon}} + 6L\varepsilon,$$

*where $c$ is an absolute constant, and*

$$Q_1(p) := \begin{cases} 1, & 1 \leq p \leq 2, \\ \max_{\substack{\lambda, \nu \in \Lambda_\varepsilon \\ \lambda \neq \nu}} \left[ \frac{\|\sum_{i=1}^N (\lambda_i - \nu_i) f_i\|_{2p-2}}{\|\sum_{i=1}^N (\lambda_i - \nu_i) f_i\|_p} \right]^{p-1}, & 2 < p < \infty. \end{cases}$$

The proof is identical to that of Corollary 1; it suffices to note only that $n_\varepsilon = \text{card}(\Lambda_\varepsilon) = (c'\varepsilon^{-1})^N$, where $c'$ is an absolute constant.

It is well known [Tsybakov (2003)] that in the problem of convex aggregation with $p = 2$ and $N \leq \varepsilon^{-1}$ the optimal (in a minimax sense) order of the remainder term is $\varepsilon\sqrt{N}$. In this particular case, our aggregation procedure achieves the indicated bound within a logarithmic in $\varepsilon^{-1}$ factor.

4.2. *Normal means model.* Consider the normal means model

$$(34) \qquad Y = \mu + \varepsilon w, \qquad \mu \in \mathbb{R}^n, \ w \sim \mathcal{N}_n(0, \Sigma),$$

where $\mu$ is an unknown vector and $\Sigma$ is the noise correlation matrix. We want to estimate $\mu$ using the observation $Y$. The model (34) is a prototype of many different nonparametric models [see, e.g., Johnstone (1998)].

Suppose that we are given a family $\Theta := \{\mu_i, i \in I_N = (1, \ldots, N)\}$ of candidate estimators of $\mu$. As before, we regard the estimators $\mu_i$, $i \in I_N$ as fixed deterministic vectors. The risk of an estimator $\hat{\mu}$ is given by $\mathbb{E}_\mu |\hat{\mu} - \mu|_p$, where $|\cdot|_p$, $p \in [1, \infty]$, stands for the standard $p$-norm in $\mathbb{R}^n$. The objective is to select a single estimator from $\Theta$ whose risk is as close as possible to that of the best estimator in $\Theta$.

The general aggregation scheme of Section 2 can be easily adapted for the outlined setup. Let $\Psi$ denote a set of probe vectors from $\mathbb{R}^n$. For $\psi \in \Psi$ define the linear functional $\ell_\mu(\psi) = \psi^T \mu$ and for every $\psi \in \Psi$ consider the following estimators of $\ell_\mu(\psi)$:

$$\hat{\ell}_\mu(\psi) = \psi^T Y, \qquad \ell_i(\psi) = \psi^T \mu_i, i \in I_N.$$

Define $\Delta_i(\psi) = \hat{\ell}_\mu(\psi) - \ell_i(\psi)$ and note that $\Delta_i(\psi) = \psi^T(\mu - \mu_i) + \varepsilon Z(\psi)$ where $Z(\psi) = \psi^T w$ is a zero-mean normal random variable with variance $|\psi|_\Sigma^2 := \psi^T \Sigma \psi$.

The aggregation procedure is defined as follows. Let $\delta \in (0,1)$, and let

$$(35) \qquad \varkappa = \varkappa(\delta, \Psi) := \min\left\{\varkappa > 0 \middle| \mathbb{P}\left(\max_{\psi \in \Psi} \frac{|Z(\psi)|}{|\psi|_\Sigma} \geq \varkappa\right) \leq \delta\right\}.$$

Let, as before, $q$ and $p$ be the conjugate exponents, and define

$$(36) \qquad \hat{M}_i := \max_{\psi \in \Psi}\left\{\frac{1}{|\psi|_q}(|\Delta_i(\psi)| - \varkappa \varepsilon |\psi|_\Sigma)\right\},$$

$$(37) \qquad \hat{i} := \arg\min_{i \in I_N} \hat{M}_i, \qquad \hat{\mu} := \mu_{\hat{i}}.$$

According to Section 2, the set of probe vectors $\Psi$ should have some "good" *norm approximation* properties. In the context of the normal means model this requirement is formulated as follows.

DEFINITION 2. Let

$$\mathcal{G} := \{g \in \mathbb{R}^n : g = \mu_i - \mu_j, i \neq j, i, j \in I_N\},$$

and let $\gamma \geq 0$. We say that the set of vectors $\Psi$ from $\mathbb{R}^n$ is $(\gamma, p)$-good if for every vector $g \in \mathcal{G}$ there is a vector $\psi_g \in \Psi$ such that

$$|\psi_g^T g - |g|_p| \leq \gamma.$$

As before we will use $(\gamma, p)$-good sets $\Psi$ in the form

$$\Psi = \{\psi | \psi = \psi_{ij} := \psi_{\mu_i - \mu_j}, i \neq j, i, j \in I_N\},$$

where $\psi_{ij}$ is a vector such that

$$|\psi_{ij}^T(\mu_i - \mu_j) - |\mu_i - \mu_j|_p| \leq \gamma.$$

Now we are in a position to establish an oracle inequality for the aggregation rule (36)–(37).

THEOREM 6. *Let $p \in [1, \infty]$, $\Psi$ be a $(\gamma, p)$-good set, $\delta \in (0, 1)$, and let $\varkappa$ be defined in (35). Assume that*

$$\max\{|\mu|_p, |\mu_1|_p, \ldots, |\mu_N|_p\} =: L < \infty.$$

*Define $i_* = \arg\min_i |\mu_i - \mu|_p$, and*

$$\Psi_* := \{\psi \in \Psi | \psi = \psi_{i_* j} = \psi_{\mu_{i_*} - \mu_j}, j \neq i_*, j \in I_N\}.$$

*Then for $\hat{\mu}$ given by (36)–(37) one has*

(38)
$$\mathbb{E}_\mu |\hat{\mu} - \mu|_p \leq \left(2 \max_{\psi \in \Psi_*} |\psi|_q + 1\right) \min_i |\mu_i - \mu|_p$$
$$+ 2\varkappa\varepsilon \max_{\psi \in \Psi_*} |\psi|_\Sigma + \gamma + 2L\delta.$$

The proof of Theorem 6 is identical to that of Theorem 1, and it is omitted.

The oracle inequality of Theorem 6 is easily specialized for specific sets of $(\gamma, p)$-good probe vectors. For example, let $p \in [1, \infty)$ and define $\widetilde{\psi}_{ij} \in \mathbb{R}^n$ by

$$\widetilde{\psi}_{ij}(k) := \frac{|\mu_i(k) - \mu_j(k)|^{p-1}}{|\mu_i - \mu_j|_p^{p-1}} \operatorname{sign}\{\mu_i(k) - \mu_j(k)\}, \qquad i, j \in I_N,$$

where $a(k)$, $k = 1, \ldots, n$, denotes the $k$th component of a generic vector $a \in \mathbb{R}^n$. Then the set of probe vectors $\widetilde{\Psi} := \{\widetilde{\psi}_{ij}, i \neq j, i, j \in I_N\}$ is $(0, p)$-good. Note also that $\widetilde{\Psi} \subset \{\psi : |\psi|_q = 1\}$.

The next result is an immediate consequence of Theorem 6.

COROLLARY 5. *Let $p \in [1, \infty)$, $\Psi = \widetilde{\Psi}$, and assume that $\Sigma$ is the identity matrix. Let $\delta = \varepsilon$; then*

$$\mathbb{E}_\mu |\hat{\mu} - \mu|_p \leq 3 \min_{i \in I_N} |\mu - \mu_i|_p + 2Q(p)\varepsilon \sqrt{2 \ln \frac{N^2}{\varepsilon}} + 2L\varepsilon,$$

*where*

$$Q(p) := \begin{cases} 1, & 2 \leq p < \infty, \\ \max_{\substack{i \in I_N \\ i \neq i_*}} \left[\frac{|\mu_{i_*} - \mu_i|_{2p-2}}{|\mu_{i_*} - \mu_i|_p}\right]^{p-1}, & 1 < p \leq 2, \\ \max_{\substack{i \in I_N \\ i \neq i_*}} [\operatorname{card}\{k : \mu_i(k) \neq \mu_{i_*}(k)\}]^{1/2}, & p = 1. \end{cases}$$

Corollary 5 shows that if $p \in [2, \infty)$, then the risk of the selected estimator is within factor 3 of the best possible risk whereas the remainder term is of the order $\varepsilon\sqrt{\ln(N^2/\varepsilon)}$. If $p \in [1, 2)$, then the remainder terms in the oracle inequalities depend on the family of aggregated estimators. The situation here is opposite to that discussed in Section 3 because of reciprocal behavior (with respect to $p$) of $\mathbb{L}_p$-norms on $[0,1]^d$ and $p$-norms in $\mathbb{R}^n$.

The aggregation procedure (36)–(37) requires prior knowledge of the noise level $\varepsilon$ and the noise covariance matrix $\Sigma$. However, (36)–(37) can be modified in the spirit of Section 3.3. Specifically, let

$$\tilde{M}_i := \max_{\psi \in \Psi} \left\{ \frac{1}{|\psi|_q} |\Delta_i(\psi)| \right\}, \tag{39}$$

$$\tilde{i} := \arg\min_{i \in I_N} \tilde{M}_i, \qquad \tilde{\mu} := \mu_{\tilde{i}}. \tag{40}$$

The next result establishes an upper bound on the accuracy of $\tilde{\mu}$.

THEOREM 7. *Let conditions of Theorem 6 hold. Then for the estimator $\tilde{\mu}$ one has*

$$\mathbb{E}_\mu |\tilde{\mu} - \mu|_p \leq \left( 2 \max_{\psi \in \Psi_*} |\psi|_q + 1 \right) \min_i |\mu_i - \mu|_p$$
$$+ 2\varkappa\varepsilon \max_{\psi \in \Psi_*} |\psi|_q \max_{\psi \in \Psi} \frac{|\psi|_\Sigma}{|\psi|_q} + \gamma + 2L\delta. \tag{41}$$

The proof is identical to that of Theorem 3 and it is omitted.

Even though the right-hand side of (41) is greater than or equal to the right-hand side of (38), $\tilde{\mu}$ can be advantageous in comparison with $\hat{\mu}$. For instance, if $p = 2$, and if the ratio of the norms $|\cdot|_\Sigma$ and $|\cdot|_2$ does not depend on $N$, then the second terms on the right-hand sides of (41) and (38) are of the same order. In this case it is advantageous to use the estimator $\tilde{\mu}$ because it does not require knowledge of $\varepsilon$ and $\Sigma$.

4.3. *Some numerical results.* A small simulation study was carried out in order to illustrate usefulness and practical potential of the proposed scheme. We investigate performance of our procedure for estimating a normal mean vector under the following two different scenarios:

(i) the vector has $K$ randomly located nonzero coefficients;
(ii) the vector has $K$ first nonzero components.

Under the first scenario thresholding estimators with properly chosen threshold will presumably perform well. In this context our selection rule provides an estimator that adapts to unknown sparsity. Recently the topic of adaptive estimation of sparse vectors has attracted much attention in the literature; we refer to Abramovich et al. (2006), Golubev (2002) and Johnstone and Silverman (2004) where further references can be found. In the second scenario projection estimators are appropriate. As we will see below, our estimator mimics the best estimator closely in both cases.

TABLE 1
*The $\mathbb{L}_2$-risk averaged over 100 replications in estimating* (i) *a normal mean vector with $K$ randomly located nonzero coefficients;* (ii) *a normal mean vector with $K$ first nonzero coefficients*

|  | $K$ | Oracle | Aggregation | Best projection estimator | Best thresholding estimator | $\widehat{K}$ |
|---|---|---|---|---|---|---|
| (i) | 5 | 2.498 | 2.726 | 4.593 | 2.499 | 5.26 |
|  | 50 | 6.446 | 6.557 | 13.994 | 6.446 | 50.08 |
|  | 250 | 11.388 | 11.559 | 19.949 | 11.388 | 292 |
|  | 500 | 13.649 | 14.378 | 24.471 | 13.649 | 613.03 |
| (ii) | 10 | 1.551 | 2.340 | 1.556 | 2.582 | 11.29 |
|  | 50 | 3.546 | 3.916 | 3.546 | 5.589 | 44.89 |
|  | 250 | 8.608 | 8.955 | 8.608 | 11.337 | 283.69 |
|  | 500 | 11.200 | 11.230 | 11.200 | 14.566 | 497.33 |

Conditions of our numerical experiments are as follows. We consider the normal means model (34) with $n = 1000$ and $\Sigma$ being the identity matrix. In the first scenario components of the unknown vector $\mu$ are assumed to be zero except $K = 5, 50, 250, 500$ randomly chosen locations where they take a specified value $m = 2$. In the second scenario the unknown vector $\mu$ has first $K = 10, 50, 250, 500$ nonzero components that are generated as independent standard normal random variables. In both scenarios the results are averaged over 100 replications for each value of $K$.

In our experiments we use two samples (random vectors) $Y_1$ and $Y_2$: the first one $Y_1 \sim \mathcal{N}_{1000}(\mu, \varepsilon_1^2 I)$, $\varepsilon_1 = 0.5$, is used for construction of estimators, while the second one $Y_2 \sim \mathcal{N}_{1000}(\mu, \varepsilon_2^2 I)$, $\varepsilon_2 = 1$, is for the aggregation purposes. The collection $\Theta$ contains 20 estimators $\hat{\mu}_1, \ldots, \hat{\mu}_{20}$:

- 10 projection estimators $\hat{\mu}_i$, $i = 1, \ldots, 10$,

$$\hat{\mu}_i(k) = Y_1(k) 1(k \leq \text{ord}_i), \qquad k = 1, \ldots, 1000,$$

with $\text{ord} = (5, 10, 20, 50, 100, 200, 300, 500, 700, 800)$.
- 10 thresholding estimators $\hat{\mu}_i$, $i = 11, \ldots, 20$,

$$\hat{\mu}_i(k) = Y_1(k) 1\{|Y_1(k)| \geq \varepsilon_1 \sqrt{2 \ln(n/t_{i-10})}\}, \qquad k = 1, \ldots, 1000,$$

where $t = (1, n^{1/4}, n^{1/2}, n^{3/4}, n^{5/6}, n^{7/8}, n^{9/10}, n^{15/16}, n^{31/32}, n^{63/64})$.

The estimators are aggregated on the basis of the second sample $Y_2$ using the modified procedure (39)–(40) with $p = 2$.

Table 1 reports on the average $\mathbb{L}_2$-risk of the proposed aggregation procedure (Aggregation), and the average $\mathbb{L}_2$-risks of three oracles that know the vector to be estimated and select: (a) the best estimator (`Oracle`) in the collection; (b) the best projection estimator in the collection; and (c) the best thresholding estimator in the collection. The last column $\widehat{K}$ displays the average number of nonzero coefficients in the selected estimate. Part (i) of the table presents results for the first scenario while part (ii) corresponds to the second scenario.

The results indicate that in estimating sparse vectors [part (i) of the table] in almost all replications thresholding estimators outperform the projection estimators. The situation is opposite for vectors with nonzero first coefficients [part (ii) of the table]: here projection estimators perform better. In both cases our aggregation procedure follows closely the best estimator from the collection for all values of $K$. The results in the last column also suggest that the aggregation procedure recovers a sparsity pattern of the estimated vector.

Additional insight into performance of the aggregation procedure is gained from Figures 1 and 2. These figures show typical behavior of the procedure under scenarios (i) and (ii). The rows (a)–(d) of the diagrams in Figures 1 and 2 correspond to different values of the parameter $K$. In each replication the competing estimators $\hat{\mu}_i$, $i = 1, \ldots, 20$, were ranked according to their performance measured by the $\mathbb{L}_2$-risk. The barplots in the left column of the figures display the number of replications out of 100 where the aggregation procedure selects the estimator with ranks $1, 2, \ldots, 20$. The diagrams in the middle column of Figures 1 and 2 show how many times the estimators $\hat{\mu}_i$ were selected. The right column displays the $\mathbb{L}_2$-risk of all estimators averaged over 100 replications.

It is seen from the barplots in the left column of Figure 1 that in the cases $K = 5, 50, 250$ the procedure selects the best estimator in more than 65% of replications. In particular, for $K = 5$ the middle panel in the row (a) demonstrates that most of the time the procedure selects the estimators $\hat{\mu}_{11}$ and $\hat{\mu}_{12}$ (the thresholding estimators with $t = 1$ and $t = n^{1/4}$, resp.). The corresponding barplot in the right column shows that the average $\mathbb{L}_2$-risks of these two estimators are significantly smaller than those of the other estimators. Similar patterns are observed when $K$ equals 50 and 250 [the rows (b) and (c) of Figure 1]. On the other hand, in the case $K = 500$ inferior estimators are chosen more frequently. Here the procedure selects one of the seven thresholding estimators with $t \geq n^{3/4}$. As the right panel in the row (d) indicates, the average $\mathbb{L}_2$-risks of these estimators are quite close. This fact explains the shape of the barplot in the corresponding left panel.

Similar conclusions can be drawn from the barplots of Figure 2. In the case $K = 10$, according to the middle panel in the row (a), the procedure selects either the projection estimators with $\text{ord} = 5, 10, 20, 50$, or the thresholding estimators with $t = 1$, $n^{1/4}$. The right panel in the row (a) shows that the average risks of these estimators are quite close. On the other hand, when $K = 500$ [the row (d) of Figure 2], the projection estimator of the order ord $= 500$ is selected in all replications, and its average risk is significantly smaller than the risks of all other estimators.

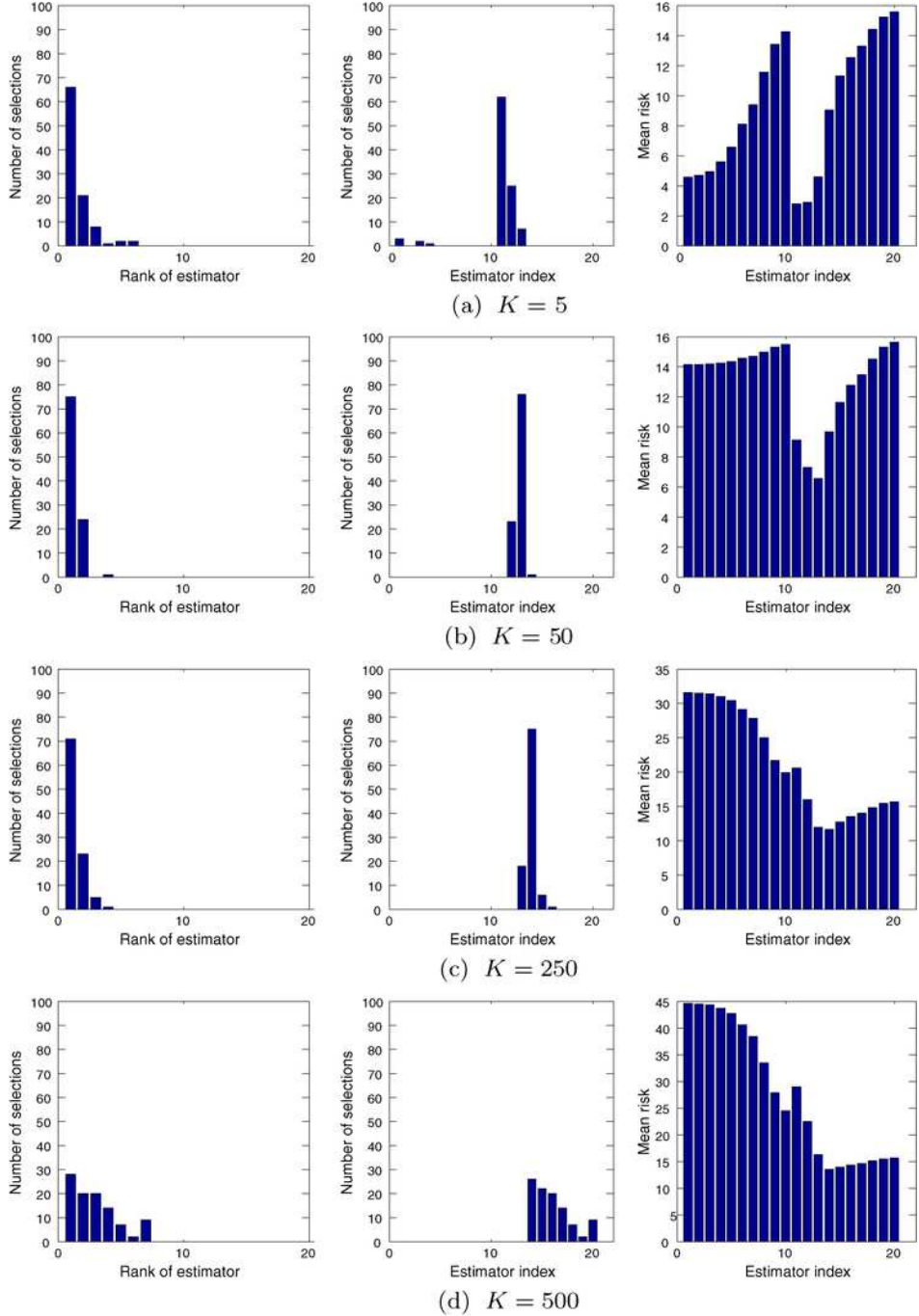

FIG. 1. *Scenario* (i). *Left column: the number of replications out of* 100 *where the procedure selects the estimator with rank* $1, 2, \ldots, 20$. *Middle column: the number of selections versus the estimator index. Right column: the average* $\mathbb{L}_2$-*risk versus the estimator index. Sparsity parameter* $K$: (a) $K = 5$; (b) $K = 50$; (c) $K = 250$; (d) $K = 500$.

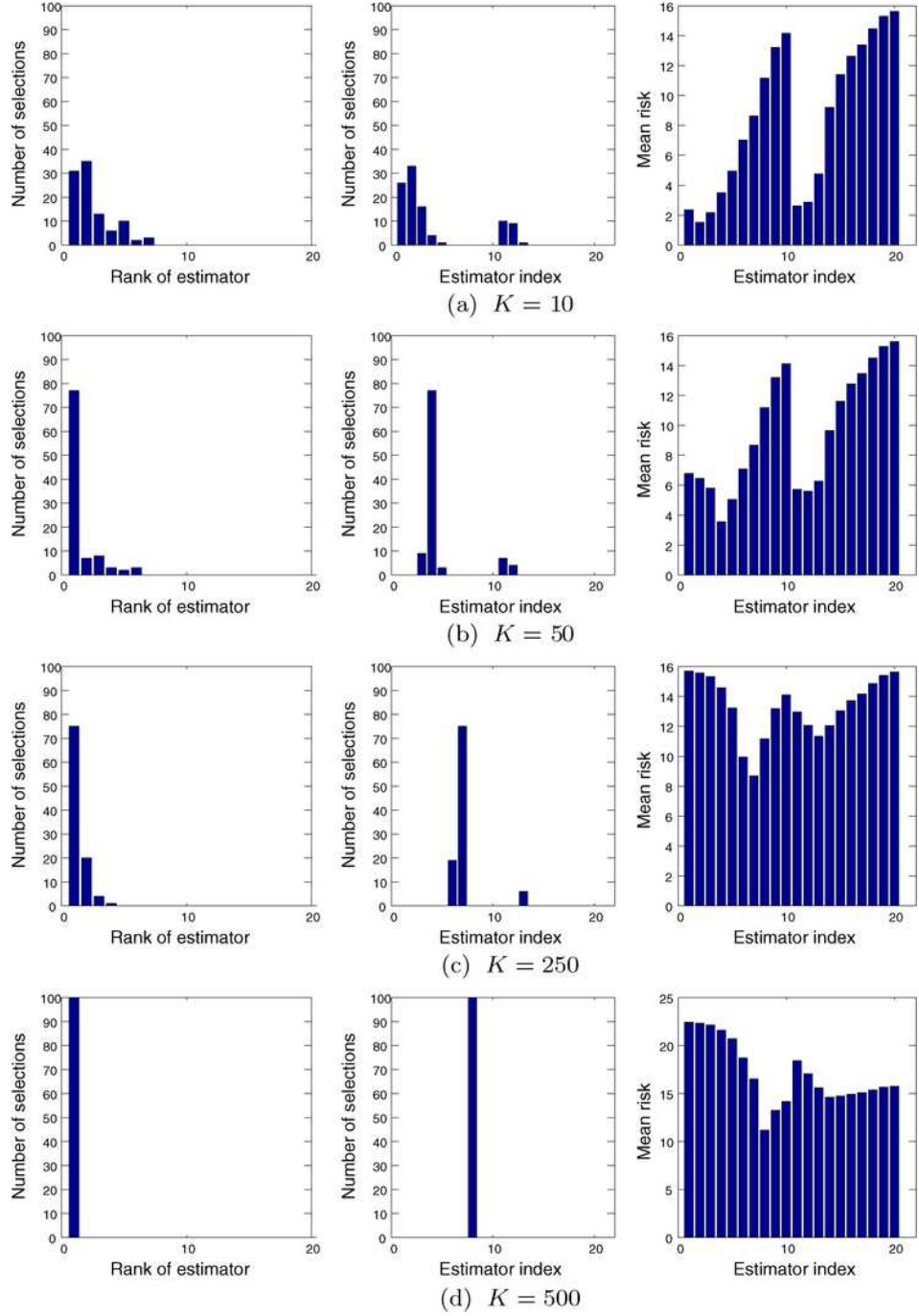

Fig. 2. *Scenario* (ii). *Left column: the number of replications out of* 100 *where the procedure selects the estimator with rank* 1, 2, ..., 20. *Middle column: the number of selections versus the estimator index. Right column: the average* $\mathbb{L}_2$-*risk versus the estimator index. The parameter* $K$: (a) $K = 5$; (b) $K = 50$; (c) $K = 250$; (d) $K = 500$.

Summing up, the shapes of the diagrams in Figures 1 and 2 and our numerical experience suggest that performance of the procedure is essentially determined by the risks of the estimators to be aggregated and by the noise level $\varepsilon_2$ at the aggregation stage. The procedure succeeds to detect the best estimator in a majority of replications when its performance is "significantly" better than the performance of the other estimators in the collection. Significance here is relative to the noise level $\varepsilon_2$ at the aggregation stage. On the other hand, if there is a large number of good estimators that perform almost equally well, the procedure makes more errors in the estimator selection. However, this does not lead to a significant increase in the risk. Our numerical experience shows also that behavior of the proposed aggregation procedure is quite reasonable for the $\mathbb{L}_1$-losses as well.

## 5. Proofs.

### 5.1. *Proofs of Theorem 1 and Corollary 1.*

PROOF OF THEOREM 1. (1) We begin with the following simple observation. Let

$$(42) \qquad A_\varkappa := \left\{\omega: \max_{\psi \in \Psi_{I_N}} \frac{|Z(\psi)|}{\|\psi\|_2} \leq \varkappa\right\},$$

where $\varkappa = \varkappa(\delta, \Psi_{I_N})$ is defined in (5). It follows from (13) and definition of $A_\varkappa$ that for any $\psi \in \Psi_{I_N}$ and $i \in I_N$

$$(43) \qquad |\Delta_i(\psi)|1(A_\varkappa) \leq \left|\int \psi(t)[f(t) - f_i(t)]\,dt\right| + \varepsilon\varkappa\|\psi\|_2.$$

Therefore

$$(44) \qquad \begin{aligned}\hat{M}_i 1(A_\varkappa) &= \max_{\psi \in \Psi_{I_N}} \frac{1}{\|\psi\|_q}[|\Delta_i(\psi)| - \varepsilon\varkappa\|\psi\|_2]1(A_\varkappa)\\ &\leq \|f - f_i\|_p \qquad \forall i \in I_N.\end{aligned}$$

(2) Write

$$\|\hat{f} - f\|_p = \|\hat{f} - f\|_p 1(A_\varkappa) + \|\hat{f} - f\|_p 1(A_\varkappa^c).$$

By definition $\mathbb{P}(A_\varkappa) \geq 1 - \delta$. Let $i_* = \arg\min_{i \in I_N} \|f - f_i\|_p$ and $f_* = f_{i_*}$; then

$$(45) \qquad \|\hat{f} - f\|_p 1(A_\varkappa) \leq \|f_* - f\|_p 1(A_\varkappa) + \|f_{i_*} - f_{\hat{i}}\|_p 1(A_\varkappa).$$

Our current goal is to bound the second term on the right-hand side of (45).

First we note that

$$\Delta_i(\psi) - \Delta_j(\psi) = \ell_{f_j}(\psi) - \ell_{f_i}(\psi)$$
$$(46) \qquad = \int \psi(t)[f_j(t) - f_i(t)]\,dt \qquad \forall i,j \in I_N, \psi \in \Psi_{I_N}.$$

By the premise of the theorem $\Psi_{I_N}$ is $(\gamma, p)$-good w.r.t. $\mathcal{G}_{I_N}$; hence there exists a probe function, say, $\psi_{i_*\hat{i}} := \psi_{f_{i_*} - f_{\hat{i}}} \in \Psi_{I_N}$ such that

$$(47) \qquad \|f_{i_*} - f_{\hat{i}}\|_p \leq \left|\int \psi_{i_*\hat{i}}(t)[f_{i_*}(t) - f_{\hat{i}}(t)]\,dt\right| + \gamma.$$

Therefore we have on the set $A_\varkappa$

$$(48) \quad \begin{aligned}
\|f_{i_*} - f_{\hat{i}}\|_p &\overset{(a)}{\leq} |\Delta_{i_*}(\psi_{i_*\hat{i}}) - \Delta_{\hat{i}}(\psi_{i_*\hat{i}})| + \gamma \\
&\leq [|\Delta_{i_*}(\psi_{i_*\hat{i}})| - \varepsilon\varkappa\|\psi_{i_*\hat{i}}\|_2] + [|\Delta_{\hat{i}}(\psi_{i_*\hat{i}})| - \varepsilon\varkappa\|\psi_{i_*\hat{i}}\|_2] \\
&\quad + 2\varepsilon\varkappa\|\psi_{i_*\hat{i}}\|_2 + \gamma \\
&\overset{(b)}{\leq} (\hat{M}_{i_*} + \hat{M}_{\hat{i}})\max_{\psi \in \Psi^*_{I_N}} \|\psi\|_q + 2\varepsilon\varkappa \max_{\psi \in \Psi^*_{I_N}} \|\psi\|_2 + \gamma \\
&\overset{(c)}{\leq} 2\hat{M}_{i_*} \max_{\psi \in \Psi^*_{I_N}} \|\psi\|_q + 2\varepsilon\varkappa \max_{\psi \in \Psi^*_{I_N}} \|\psi\|_2 + \gamma \\
&\overset{(d)}{\leq} 2\left[\max_{\psi \in \Psi^*_{I_N}} \|\psi\|_q\right]\|f - f_{i_*}\|_p + 2\varepsilon\varkappa \max_{\psi \in \Psi^*_{I_N}} \|\psi\|_2 + \gamma,
\end{aligned}$$

where (a) follows from (46) and (47), (b) is by definition of $\hat{M}_i$ and because $\psi_{i_*\hat{i}} \in \Psi^*_{I_N}$ [see (16)], (c) follows from (15) and (d) is by (44).

(3) On the set $A^c_\varkappa$ we have

$$\|\hat{f} - f\|_p 1(A^c_\varkappa) \leq \left[\|f\|_p + \max_{i \in I_N} \|f_i\|_p\right] 1(A^c_\varkappa).$$

Combining this inequality with (48) and (45) we complete the proof. □

PROOF OF COROLLARY 1. By Example 1, $\widetilde{\Psi}_{I_N}$ is $(0, p)$-good so that $\gamma = 0$ in (17). Moreover, $\|\psi\|_q = 1$ for all $\psi \in \widetilde{\Psi}_{I_N}$. Since the cardinality of $\widetilde{\Psi}_{I_N}$ equals $N(N-1)$ we have

$$\mathbb{P}\left\{\max_{\psi \in \widetilde{\Psi}_{I_N}} \frac{|Z(\psi)|}{\|\psi\|_2} \geq \varkappa\right\} \leq N^2 \exp\{-\varkappa^2/2\}.$$

It follows from the definition of $\varkappa$ and the preceding inequality that $N^2 e^{-\varkappa^2/2} \geq \delta$ so that $\varkappa \leq \sqrt{2\ln(N^2/\delta)} = \sqrt{2\ln(N^2/\varepsilon)}$.

If $p \in [1,2]$, then
$$\max_{\psi \in \widetilde{\Psi}_{I_N}} \|\psi\|_2 \leq \max_{\psi \in \widetilde{\Psi}_{I_N}} \|\psi\|_q = 1.$$

On the other hand, if $2 < p < \infty$, then in view of (11)
$$\max_{\psi \in \widetilde{\Psi}^*_{I_N}} \|\psi\|_2 = \max_{\substack{i \in I_N \\ i \neq i_*}} \left[ \frac{\|f_{i_*} - f_i\|_{2p-2}}{\|f_{i_*} - f_i\|_p} \right]^{p-1}.$$

Combining these inequalities with the statement of Theorem 1 we come to (19). □

5.2. *Proof of Theorem 2.* Let $B_i$, $i = 1,\ldots,N$ be disjoint subsets of $\mathcal{D}_0$ such that $\mathrm{mes}(B_i) = h$, $\forall i$, where $0 < h \leq 1/N$, is a given number. Here $\mathrm{mes}(\cdot)$ stands for the Lebesgue measure in $\mathbb{R}^d$. Define $\bar{f}_i(x) = L 1_{B_i}(x)$, $i \in I_N$, and $\bar{\mathcal{F}}_{I_N} = \{\bar{f}_i, i \in I_N\}$. Note that $\max_{i \in I_N} \|\bar{f}_i\|_p \leq L$ for all $p \in (2, \infty]$. If $f \in \bar{\mathcal{F}}_{I_N}$, then $\min_{i \in I_N} \|f - \bar{f}_i\|_p = \|\bar{f}_{i_*} - f\|_p = 0$. Moreover
$$\|\bar{f}_i - \bar{f}_j\|_p = (2h)^{1/p} L =: s \qquad \forall i, j \in I_N, i \neq j,$$
and
$$(49) \quad Q_1(\bar{\mathcal{F}}_{I_N}, p) = \max_{\substack{i \in I_N, \\ i \neq i_*}} \frac{\|\bar{f}_{i_*} - \bar{f}_i\|_{2p-2}^{p-1}}{\|\bar{f}_{i_*} - \bar{f}_i\|_p^{p-1}} = (2h)^{1/p - 1/2},$$
$$Q_3(\bar{\mathcal{F}}_{I_N}, p) = \max_{\substack{i \in I_N \\ i \neq i_*}} \frac{\|\bar{f}_{i_*} - \bar{f}_i\|_p}{\|\bar{f}_{i_*} - \bar{f}_i\|_2} = (2h)^{1/p - 1/2}.$$

It is immediately seen that for a chosen family of functions one has
$$Q_4(\bar{\mathcal{F}}_{I_N}, \gamma) = \frac{\gamma (2h)^{1/2}}{\gamma (2h)} = (2h)^{-1/2} \qquad \forall \gamma > 0,$$

which coincides with (49) for $p = \infty$. Denote $K_p := (2h)^{1/p - 1/2}$, $p \in (2, \infty]$.

Let $\tilde{f} : \mathcal{Y}_\varepsilon \to \bar{\mathcal{F}}_{I_N}$ be an arbitrary selection rule. We have

$$(50) \quad \sup_{f \in \bar{\mathcal{F}}_{I_N}} \mathbb{E}_f \|\tilde{f} - f\|_p \geq \frac{s}{2} \max_{i \in I_N} \mathbb{P}_i \left\{ \|\tilde{f} - \bar{f}_i\|_p \geq \frac{s}{2} \right\} \geq \frac{s}{2} \max_{i \in I_N} \mathbb{P}_i \{\tilde{i} \neq i\},$$

where $\mathbb{P}_i = \mathbb{P}_{\bar{f}_i}$ probability measure of the observation $\mathcal{Y}_\varepsilon$ associated with $f = \bar{f}_i$, and $\tilde{i} : \mathcal{Y}_\varepsilon \to \{1,\ldots,N\}$ is the decision rule that selects function $\bar{f}_i$ closest to $\tilde{f}$ in the $\mathbb{L}_p$-norm.

Let $K(\mathbb{P}_i, \mathbb{P}_j)$ denote the Kullback–Leibler divergence between $\mathbb{P}_i$ and $\mathbb{P}_j$:
$$K(\mathbb{P}_i, \mathbb{P}_j) = \frac{1}{2\varepsilon^2} \|\bar{f}_i - \bar{f}_j\|_2^2 = \frac{hL^2}{\varepsilon^2} \qquad \forall i, j \in I_N, i \neq j.$$

Then by the Fano inequality [see, e.g., Devroye (1987), Section 5.9]

$$\max_{i \in I_N} \mathbb{P}_i\{\tilde{i} \neq i\} \geq 1 - \frac{hL^2\varepsilon^{-2} + \ln 2}{\ln(N-1)}.$$

Choosing

$$h = h_* = \frac{\varepsilon^2}{L^2}\left(\frac{5}{6}\ln(N-1) - \ln 2\right) \geq \frac{\varepsilon^2}{6L^2}\ln(N-1)$$

(the last inequality follows from $N > 3$), we obtain that $\max_i \mathbb{P}_i\{\tilde{i} \neq i\} \geq 1/6$. Note that condition $\varepsilon \leq L(N \ln N)^{-1/2}$ implies $h_* \leq 1/N$ so that the sets $B_i$ are indeed disjoint, as assumed. Hence (50) yields

$$\sup_{f \in \bar{\mathcal{F}}_{I_N}} \mathbb{E}_f \|\tilde{f} - f\|_p \geq \frac{L}{12}(2h_*)^{1/p}$$

$$= \frac{K_p}{12}L(2h_*)^{1/2} \geq \frac{K_p}{12\sqrt{3}}\varepsilon\sqrt{\ln(N-1)}.$$

This completes the proof.

5.3. *Proof of Theorem 3.* The proof goes along the same lines as the proof of Theorem 1; below we indicate only the differences. We use the same notation as in the proof of Theorem 1.

First we note that for all $i \in I_N$

$$\tilde{M}_i 1(A_\varkappa) = \max_{\psi \in \Psi_{I_N}} \left\{\frac{1}{\|\psi\|_q}|\Delta_i(\psi)|\right\} 1(A_\varkappa) \leq \|f - f_i\|_p + \varepsilon\varkappa \max_{\psi \in \Psi_{I_N}} \frac{\|\psi\|_2}{\|\psi\|_q}.$$

Because $\Psi_{I_N}$ is $(\gamma, p)$-good, there is a probe function, say, $\psi_{i_*\tilde{i}} \in \Psi_{I_N}$ such that

$$\|f_{i_*} - f_{\tilde{i}}\|_p \leq \left|\int \psi_{i_*\tilde{i}}(t)[f_{i_*}(t) - f_{\tilde{i}}(t)]\,dt\right| + \gamma.$$

Then, similarly to (48), we have on the set $A_\varkappa$

$$\|f_{i_*} - f_{\tilde{i}}\|_p \leq |\Delta_{i_*}(\psi_{i_*\tilde{i}}) - \Delta_{\tilde{i}}(\psi_{i_*\tilde{i}})| + \gamma$$

$$\leq \|\psi_{i_*\tilde{i}}\|_q(\tilde{M}_{i_*} + \tilde{M}_{\tilde{i}}) + \gamma$$

$$\leq 2\tilde{M}_{i_*} \max_{\psi \in \Psi_{I_N}^*} \|\psi\|_q + \gamma$$

$$\leq 2\max_{\psi \in \Psi_{I_N}^*} \|\psi\|_q \left(\|f - f_{i_*}\|_p + \varepsilon\varkappa \max_{\psi \in \Psi_{I_N}} \frac{\|\psi\|_2}{\|\psi\|_q}\right) + \gamma.$$

This leads to the inequality (28).

5.4. *Proof of Theorem 4.* Throughout the proof $\langle \cdot, \cdot \rangle$ denotes the standard inner product in $\mathbb{L}_2(\mathcal{D}_0)$.

We start with the following simple observation. Let $f_{i_*}$ be the best estimator in the family $\mathcal{F}_{I_N}$, that is, $i_* = \arg\min_{i \in I_N} \|f_i - f\|_2$. Since for any $j \in I_N$
$$\|f_{i_*} - f\|_2^2 = \|f_j - f\|_2^2 + \|f_{i_*} - f_j\|_2^2 + 2\langle f_{i_*} - f_j, f_j - f \rangle$$
and $\|f_{i_*} - f\|_2 \leq \|f_j - f\|_2$, then
$$\|f_{i_*} - f_j\|_2^2 + 2\langle f_{i_*} - f_j, f_j - f \rangle = 2\langle f_{i_*} - f_j, \tfrac{1}{2}(f_{i_*} + f_j) - f \rangle$$
$$= 2\langle f_{i_*} - f_j, u_{i_*j} - f \rangle \leq 0 \qquad \forall j \in I_N,$$
or, equivalently,
$$\max_{j \in I_N} \langle \psi_{i_*j}, u_{i_*j} - f \rangle \leq 0. \tag{51}$$

We have
$$\|\bar{f} - f\|_2^2 = \|f_{i_*} - f\|_2^2 + 2\langle f_{\bar{i}} - f_{i_*}, \tfrac{1}{2}(f_{\bar{i}} + f_{i_*}) - f \rangle$$
$$\stackrel{(a)}{=} \|f_{i_*} - f\|_2^2 + 2\|f_{\bar{i}} - f_{i_*}\|_2 \langle \psi_{\bar{i}i_*}, u_{\bar{i}i_*} - f \rangle$$
$$= \|f_{i_*} - f\|_2^2 + 2\|f_{\bar{i}} - f_{i_*}\|_2 \left\{ \langle \psi_{\bar{i}i_*}, u_{\bar{i}i_*} \rangle - \int \psi_{\bar{i}i_*}(t) Y_\varepsilon(dt) \right\}$$
$$+ 2\|f_{\bar{i}} - f_{i_*}\|_2 \varepsilon Z(\psi_{\bar{i}i_*}) \tag{52}$$
$$\stackrel{(b)}{\leq} \|f_{i_*} - f\|_2^2 + 2\|f_{\bar{i}} - f_{i_*}\|_2 \bar{M}_{\bar{i}} + 2\|f_{\bar{i}} - f_{i_*}\|_2 \varepsilon Z(\psi_{\bar{i}i_*})$$
$$\stackrel{(c)}{\leq} \|f_{i_*} - f\|_2^2 + 2\|f_{\bar{i}} - f_{i_*}\|_2 \bar{M}_{i_*} + 2\|f_{\bar{i}} - f_{i_*}\|_2 \varepsilon Z(\psi_{\bar{i}i_*}),$$
where $Z(\psi) = \int \psi(t) W(dt)$, (a) is by definition of $u_{ij}$ and $\psi_{ij}$, (b) is by definition of $\bar{M}_i$, and (c) follows from the definition of $\bar{i}$.

Now we note that
$$\bar{M}_{i_*} \leq \max_{j \in I_N} \langle \psi_{i_*j}, u_{i_*j} - f \rangle + \varepsilon \max_{j \in I_N} Z(\psi_{i_*j}) \leq \varepsilon \max_{j \in I_N} Z(\psi_{i_*j}),$$
where the last inequality is a consequence of (51). Therefore it follows from (52) and $Z(\psi_{ij}) = -Z(\psi_{ji})$, $\forall i, j$ that
$$\|f_{\bar{i}} - f\|_2^2 \leq \|f_{i_*} - f\|_2^2 + 4\|f_{\bar{i}} - f_{i_*}\|_2 \varepsilon \max_{j \in I_N} |Z(\psi_{i_*j})|.$$
Hence by the triangle inequality
$$\|f_{\bar{i}} - f\|_2^2 - \|f_{i_*} - f\|_2^2 \leq 4(\|f_{\bar{i}} - f\|_2 + \|f_{i_*} - f\|_2) \varepsilon \max_{j \in I_N} |Z(\psi_{i_*j})|$$
and finally
$$\|\bar{f} - f\|_2 \leq \|f_{i_*} - f\|_2 + 4\varepsilon \max_{j \in I_N} |Z(\psi_{i_*j})|.$$
Taking the expectation we complete the proof.

5.5. *Proofs of Lemma 1 and Theorem 5.*

PROOF OF LEMMA 1. Let $g \in \mathcal{G}_\Lambda$, that is, for some $\lambda, \nu \in \Lambda$ one has $g = \sum_{i=1}^N (\lambda_i - \nu_i) f_i$. There exist $\tilde{\lambda}, \tilde{\nu} \in \Lambda_\eta$ such that $|\tilde{\lambda} - \lambda|_1 \leq \eta$ and $|\tilde{\nu} - \nu|_1 \leq \eta$. Define $\tilde{g} = \sum_{i=1}^N (\tilde{\lambda}_i - \tilde{\nu}_i) f_i$; by definition, $\tilde{g} \in \mathcal{G}_{\Lambda_\eta}$. Because $\Psi$ is $(0, p)$-good with respect to $\mathcal{G}_{\Lambda_\eta}$, there exists $\psi = \psi_{\tilde{g}} \in \Psi$ such that

$$\int \psi_{\tilde{g}}(t) \tilde{g}(t)\, dt = \|\tilde{g}\|_p.$$

With this representer $\psi_{\tilde{g}}$ applied to $g \in \mathcal{G}_\Lambda$ we obtain

$$\int \psi_{\tilde{g}}(t) g(t)\, dt = \|\tilde{g}\|_p + \int \psi_{\tilde{g}}(t)[g(t) - \tilde{g}(t)]\, dt,$$

and therefore

$$\left| \int \psi_{\tilde{g}}(t) g(t)\, dt - \|g\|_p \right| \leq |\|\tilde{g}\|_p - \|g\|_p| + \left| \int \psi_{\tilde{g}}(t)[g(t) - \tilde{g}(t)]\, dt \right|$$

$$\leq \|\tilde{g} - g\|_p + \|\psi_{\tilde{g}}\|_q \|\tilde{g} - g\|_p$$

$$= (1 + \|\psi_{\tilde{g}}\|_q) \|\tilde{g} - g\|_p.$$

To complete the proof it is sufficient to note that

$$\tilde{g}(t) - g(t) = \sum_{i=1}^N (\tilde{\lambda}_i - \lambda_i) f_i(t) - \sum_{i=1}^N (\tilde{\nu}_i - \nu_i) f_i(t);$$

hence

$$\|\tilde{g} - g\|_p \leq \sum_{i=1}^N [|\tilde{\lambda}_i - \lambda_i| + |\tilde{\nu}_i - \nu_i|] \|f_i\|_p \leq 2L\eta. \qquad \square$$

PROOF OF THEOREM 5. The proof goes along the same lines as the proof of Theorem 1; here we indicate only the main differences. First we note that similarly to (44) one has

$$\max_{\psi \in \Psi_{\Lambda_\varepsilon}} \frac{1}{\|\psi\|_q} [|\Delta_\lambda(\psi)| - \varepsilon \varkappa \|\psi\|_2] 1(A_\varkappa) \leq \|f - F_\lambda\|_p \qquad \forall \lambda \in \Lambda,$$

where $A_\varkappa$ is the event defined in (42) with $\max_{\psi \in \Psi_{I_N}}$ replaced by $\max_{\psi \in \Psi_{\Lambda_\varepsilon}}$.

Define $\lambda_* = \arg\min_\lambda \|f - F_\lambda\|_p$. The main difference with the proof of Theorem 1 is that now the set of probe functions $\Psi_{\Lambda_\varepsilon}$ is $(\gamma, p)$-good with respect to $\mathcal{G}_\Lambda$ with $\gamma$ given by (32), and the inequality (47) holds for some representer, say $\psi_{\hat{\lambda}, \nu}$, with $\nu \in \Lambda_\varepsilon$. In contrast to the proof of Theorem 1, in general $\nu \neq \lambda_*$, because $\lambda_*$ does not necessarily belong to $\Lambda_\varepsilon$. This implies that in the resulting oracle inequality we have maxima over $\psi \in \Psi_{\Lambda_\varepsilon}$, and not over the subset of $\Psi_{\Lambda_\varepsilon}$ related to $\lambda_*$. All other details of the proof remain unchanged. $\square$

**Acknowledgments.** I would like to thank A. Juditsky for useful discussions and suggestions, and an anonymous referee for comments that prompted me to improve the presentation of the numerical results.

DEPARTMENT OF STATISTICS
UNIVERSITY OF HAIFA
31905 HAIFA
ISRAEL
E-MAIL: goldensh@stat.haifa.ac.il